%% file: diff_llcvx.tex
\documentclass[12pt]{article}
\usepackage{geometry}

\usepackage{amssymb,amsthm,amsmath,amsfonts,mathtools}

\usepackage{graphicx}
\usepackage{epstopdf}
\usepackage{subcaption}
\usepackage{adjustbox}
\usepackage{xcolor}
\usepackage{tikz-qtree}
\usepackage{textcomp}
\usepackage[T1]{fontenc}

\usepackage{listings}
\definecolor{seagreen}{rgb}{0.18, 0.55, 0.34}
\definecolor{mediumviolet-red}{rgb}{0.78, 0.08, 0.52}
\definecolor{khaki}{rgb}{0.94, 0.9, 0.55}

\input{defs.tex}

\usepackage[
    natbib=true,
    citestyle=alphabetic,
    style=alphabetic,
    maxnames=3,
    minnames=3,
    maxbibnames=99,
    backend=bibtex,
    urldate=iso8601,
    firstinits=true,
    isbn=false,
    doi=false,
    date=year,
    sorting=nty
]{biblatex}
\DeclareFieldFormat[article,inbook,incollection,inproceedings,patent,thesis,
  unpublished]{citetitle}{#1}
\DeclareFieldFormat[article,inbook,incollection,inproceedings,patent,thesis,
  unpublished]{title}{#1}
\renewbibmacro{in:}{%
  \ifentrytype{article}{}{\printtext{\bibstring{in}\space}}}
\renewcommand*{\labelalphaothers}{\textsuperscript{+}}

\usepackage[colorlinks,allcolors=blue,bookmarks=false,hypertexnames=true]{hyperref}
\usepackage{breakurl}
\usepackage{url}

\DeclarePairedDelimiter{\norm}{\lvert \lvert}{\rvert \rvert}

\newcommand{\Smap}{\mathcal{S}}
\newcommand{\Cmap}{\mathcal{C}}
\newcommand{\Rmap}{\mathcal{R}}
\newcommand{\DD}{\mathsf{D}}
\newcommand{\dd}{\mathsf{d}}

\makeatother

\begin{document}

\title{Differentiating through Log-Log Convex Programs}

\author{
Akshay Agrawal \\ \texttt{
\small akshayka@cs.stanford.edu} \and
Stephen Boyd \\ \texttt{\small boyd@stanford.edu}
}

\maketitle

\begin{abstract}
We show how to
efficiently compute the derivative (when it exists) of the solution map of
\emph{log-log convex programs} (LLCPs). These are nonconvex, nonsmooth optimization
problems with positive variables that become convex when the variables,
objective functions, and constraint functions are replaced with their logs. We
focus specifically on LLCPs generated by disciplined geometric
programming, a grammar consisting of a set of atomic functions with known
log-log curvature and a composition rule for combining them. We represent a
parametrized LLCP as the composition of a
smooth transformation of parameters, a convex optimization problem, and an
exponential transformation of the convex optimization problem's solution. 
The derivative of this composition can be computed efficiently,
using recently developed methods for differentiating through convex
optimization problems.
We implement our method in CVXPY, a Python-embedded modeling language and
rewriting system for convex optimization. In just a few lines of code, a user
can specify a parametrized LLCP, solve it, and evaluate the
derivative or its adjoint at a vector. This makes it possible to
conduct sensitivity analyses of solutions, given perturbations to the
parameters, and to compute the gradient of a function of the solution with
respect to the parameters. We use the adjoint of the derivative to implement
differentiable log-log convex optimization layers in PyTorch and TensorFlow.
Finally, we present applications to designing queuing systems and fitting
structured prediction models.
\end{abstract}

\section{Introduction}\label{sec-intro}

\subsection{Log-log convex programs}
A log-log convex program (LLCP) is a mathematical optimization problem
in which the variables are positive, the objective and inequality constraint
functions are \emph{log-log convex}, and the equality constraint functions
are \emph{log-log affine}. A function $f : D \subseteq \reals^n_{++} \to
\reals_{++}$ ($\reals_{++}$ denotes the positive reals) is log-log convex if
for all $x$, $y \in D$ and $\theta
\in [0, 1]$,
\[
f(x^{\theta} \circ y^{1-\theta}) \leq f(x)^\theta \circ f(y)^{1 - \theta},
\]
and $f$ is log-log affine if the inequality holds with equality
(the powers are meant elementwise and $\circ$ denotes the elementwise
product). Similarly, $f$ is log-log concave if the inequality holds when its
direction is reversed. A LLCP has the standard form
\begin{equation}\label{e-llcp}
\begin{array}{ll}
\mbox{minimize} & f_0(x) \\
\mbox{subject to} & f_i(x) \leq 1, \quad i=1, \ldots, m_1 \\
                  & g_i(x) = 1, \quad i=1, \ldots, m_2,
\end{array}
\end{equation}
where $x \in \reals^{n}_{++}$ is the variable, the
functions $f_i$ are log-log convex, and $g_i$ are log-log affine
\citep{agrawal2019dgp}. A value of the variable is a \emph{solution} of the
problem if it minimizes the objective function, among all values satisfying the
constraints.

The problem~(\ref{e-llcp}) is not convex, but it can be readily transformed to
a convex optimization problem. We make the change of variables $u = \log x$,
replace each function $f$ appearing in the LLCP with its log-log
transformation, defined by $F(u) = \log f(e^u)$, and replace the right-hand
sides of the constraints with $0$. Because the log-log transformation of a
log-log convex function is convex, we obtain an equivalent convex optimization
problem. This means that LLCPs can be solved efficiently and globally, using
standard algorithms for convex optimization \citep{boyd2004convex}.

The class of log-log convex programs is large, including many interesting
problems as special cases. Geometric programs (GPs) form a well-studied
subclass of LLCPs; these are LLCPs in which the equality constraint functions
are \emph{monomials}, of the form $x \mapsto c x_1^{a_1} x_2^{a_2} \cdots
x_n^{a_n}$, with $a_1, \ldots a_n \in \reals$ and $c \in \reals_{++}$,
and the objective and inequality constraint functions are sums of monomials,
called \emph{posynomials} \citep{duffin1967geometric, boyd2007tutorial}. GPs
have found application in digital and analog circuit design \citep{boyd2005,
hershenson2001opamp, li2004, xu2004}, aircraft design
\citep{hoburg2014geometric, brown2018vehicle, saab2018robust}, epidemiology
\citep{preciado2014, ogura2016efficient}, chemical engineering
\citep{clasen1984}, communication systems \citep{kandukuri2002,
chiang2005geometric, chiang2007power}, control \cite{ogura2019geometric},
project management \cite{ogura2019resource}, and data fitting
\citep{hoburg2016data, calafiore2019log}; for more, see
\citep[\S1.1]{agrawal2019dgp} and \citep[\S10.3]{boyd2007tutorial}

In this paper we consider LLCPs in which the objective and constraint functions
are parametrized, and we are interested in computing how a solution to an LLCP
changes with small perturbations to the parameters. For example, in a GP, the
parameters are the coefficients and exponents appearing in
monomials and posynomials. While sensitivity analysis of GPs is
well-studied \citep{dinkel1977sensitivity, dembo1982sensitivity,
kyparisis1988sensitivity, kyparisis1990sensitivity}, sensitivity analysis of
LLCPs has not, to our knowledge, previously appeared in the literature. Our
emphasis in this paper is on practical computation, instead of a theoretical
characterization of the differentiability of the solution map.

\subsection{Solution maps and sensitivity analysis}
An optimization problem can be viewed as a multivalued function mapping
parameters to the set of solutions; this set might contain zero, one, or many
elements. In neighborhoods where this \emph{solution map} is single-valued,
it is an implicit function of the parameters \citep{dontchev2009implicit}. In
these neighborhoods it is meaningful to discuss how perturbations in the
parameters affect the solution. The point of this paper is to efficiently
calculate the sensitivity of the solution of an LLCP to these perturbations, by
implicitly differentiating the solution map; this calculation also lets us
compute the gradient of a scalar-valued function of the solution, with respect
to the parameters.

There is a large body of work on the sensitivity analysis of 
optimization problems, going back multiple decades. Early papers include
\citep{fiacco1968nonlinear} and \citep{fiacco1976sensitivity}, which apply the
implicit function theorem to the first-order KKT conditions of a nonlinear
program with twice-differentiable objective and constraint functions. A
similar method was applied to GPs in \citep{kyparisis1988sensitivity,
kyparisis1990sensitivity}. Much of the work on sensitivity analysis of GPs
focuses on the special structure of the dual program (\eg,
\citep{duffin1967geometric, dembo1982sensitivity}). Various results on sensitivity
analyses of optimization problems, including nonlinear programs, semidefinite
programs, and semi-infinite programs, are collected in
\citep{bonnans2000perturbation}.

Recently, a series of papers developed methods to calculate
the derivative of convex optimization problems, in which the objective and
constraint functions may be nonsmooth. The paper \citep{busseti2019solution}
phrased a convex cone program as the problem of finding a zero of a certain
residual map, and the papers \citep{diffcp2019, amos2019differentiable} showed
how to differentiate through cone programs (when certain regularity conditions
are satisfied) by a straightforward application of the implicit function
theorem to this residual map. In \citep{agrawal2019differentiable}, a method
was developed to differentiate through high-level descriptions of convex
optimization problems, specified in a domain-specific language for convex
optimization. The method from \citep{agrawal2019differentiable} reduces
convex optimization problems to cone programs in an efficient and
differentiable way. The present paper can be understood as an analogue of
\citep{agrawal2019differentiable} for LLCPs.

\subsection{Domain-specific languages for optimization}\label{sub-s-dsl}
Log-log convex functions satisfy an important composition rule, analogous
to the composition rule for convex functions. Suppose
$h : D \subseteq \reals^{m}_{++} \to \reals_{++} \cup \{+\infty\}$ is log-log
convex, and let $[I_1, I_2, I_3]$ be a partition of $\{1, 2, \ldots, m\}$
such that $f$ is nondecreasing in the arguments index by 
$I_1$ and nonincreasing in the arguments indexed
by $I_2$. If $g$ maps a subset of $\reals^n_{++}$ into $\reals^{m}_{++}$
such that its components $g_i$  are log-log convex for $i \in I_1$, log-log
concave for $i \in I_2$, and log-log affine for $i \in I_3$,
then the composition
\[
f = h \circ g
\]
is log-log convex. An analogous rule holds for log-log concave functions.

When combined with a set of atomic functions with known log-log curvature
and per-argument monotonicities, this composition rule defines
a \emph{grammar} for log-log convex functions, \ie, a rule for combining
atomic functions to create other functions with verifiable log-log curvature.
This is the basis of disciplined geometric programming (DGP), a grammar for
LLCPs \citep{agrawal2019dgp}. In addition to compositions of atomic functions,
DGP also includes LLCPs as valid expressions, permitting the minimization of a
log-log convex function (or maximization of a log-log concave function),
subject to inequality constraints $f(x) \leq g(x)$, where $f$ is log-log
convex and $g$ is log-log concave, and equality constraints $f(x)
= g(x)$, where $f$ and $g$ are log-log affine. In \S\ref{sub-s-parameters}, we
extend DGP to include parametrized LLCPs.

The class of DGP problems is a subclass of LLCPs. Depending on the choice of
atomic functions, or \emph{atoms}, this class can be made quite large. For
example, taking powers, products, and sums as the atoms yields GP; adding the
maximum operator
yields generalized geometric programming (GGP). Several
other atoms can be added, such as the exponential function, the
logarithm, and functions of elementwise positive matrices, yielding a subclass
of LLCPs strictly larger than GGP (see \cite[\S3]{agrawal2019dgp} for
examples). In this paper we restrict our attention to LLCPs generated by DGP;
this is not a limitation in practice, since the atom library is extensible.

DGP can be used as the grammar for a domain-specific language (DSL) for log-log
convex optimization. A DSL for log-log convex optimization parses LLCPs written
in a human readable form, rewrites them into canonical forms, and compiles the
canonical forms into numerical data for low-level numerical solvers. Because
valid problems are guaranteed to be LLCPs, the DSL can guarantee that the
compilation and the numerical solve are correct (\ie, valid problems can be
solved globally). By abstracting away the numerical solver, DSLs make
optimization accessible, vastly decreasing the time between formulating
a problem and solving it with a computer. Examples of DSLs for LLCPs
include CVXPY \citep{diamond2016cvxpy, agrawal2018rewriting} and CVXR
\citep{fu2017cvxr}; additionally, CVX \citep{cvx}, GPKit
\citep{burnell2020gpkit}, and Yalmip \citep{yalmip} support GPs.

Finally, we mention that modern DSLs for convex optimization are based on
disciplined convex programming (DCP) \citep{grant2006disciplined}, which is
analogous to DGP\@. CVXPY, CVXR, CVX, and Convex.jl \citep{convexjl}
support convex optimization using DCP as the grammar. In
\cite{agrawal2019differentiable}, a method for differentiating
through parametrized DCP problems was developed; this method was implemented in
CVXPY, and PyTorch \cite{paszke2019pytorch} and TensorFlow
\cite{abadi2016tensorflow, agrawal2019tensorflow} wrappers for differentiable
CVXPY problems were implemented in a Python package called CVXPY Layers.

\subsection{This paper} In this paper we describe how to efficiently
compute the derivative of a LLCP, when it exists, specifically considering
LLCPs generated by DGP\@. In particular, we show how to evaluate the derivative
(and its adjoint) of the solution map of an LLCP at a vector.
To do this, we first extend the DGP ruleset to include parameters as atoms,
in \S\ref{s-dgp}. Then, in \S\ref{s-solution-map}, we
represent a parametrized DGP problem by the composition of a smooth
transformation of parameters, a parametrized DCP problem, and an exponential
transformation of the DCP problem's solution. We differentiate through
this composition using recently developed methods from
\citep{busseti2019solution, diffcp2019, agrawal2019differentiable} to
differentiate through the DCP problem. Unlike prior work on sensitivity
analysis of GPs, in which the objective and constraint functions are smooth,
our method extends to problems with nonsmooth objective and constraints.

We implement the derivative of LLCPs as an abstract linear operator in CVXPY\@.
In just a few lines of code, users can conduct first-order sensitivity analyses
to examine how the values of variables would change given small perturbations
of the parameters. Using the adjoint of the derivative operator, users can
compute the gradient of a function of the solution to a DGP
problem, with respect to the parameters. For convenience, we implement PyTorch
and TensorFlow wrappers of the adjoint derivative in CVXPY Layers, making it
easy to use log-log convex optimization problems as tunable layers in
differentiable programs or neural networks.
For our implementation, the overhead in differentiating through the DSL
(which rewrites the high-level description of a problem into low-level
numerical data for a solver, and retrieves a solution for the original problem
from a solution from the solver) is small compared to the time spent in
the numerical solver. In particular, the mapping from the transformed
parameters to the numerical solver data is affine and can be
represented compactly by a sparse matrix, and so can be evaluated quickly. Our implementation is described and
illustrated with usage examples in \S\ref{s-implementation}. 

In \S\ref{s-example}, we present two simple examples, in which we
apply a sensitivity analysis to the design of an $M/M/N$ queuing system,
and fit a structured prediction model.

\subsection{Related work}
\paragraph{Automatic differentiation.} Automatic differentiation (AD) is a
family of methods that use the chain rule to algorithmically compute exact
derivatives of compositions of differentiable functions, dating back to the
1950s \cite{beda1959}. There are two main types of AD\@. Reverse-mode AD
computes the gradient of a scalar-valued composition of differentiable
functions by applying the adjoint of the intermediate derivatives to the
sensitivities of their outputs, while forward-mode AD computes applies the
derivative of the composition to a vector of perturbations in the inputs
\citep{griewank2008}. In AD, derivatives are typically implemented as
\emph{abstract linear maps}, \ie, methods for applying the derivative and
its adjoint at a vector; the derivative matrices of the intermediate
functions are not materialized, \ie, formed or stored as arrays.

Recently, many high-quality open-source implementations of AD were made
available. Examples include PyTorch \citep{paszke2019pytorch}, TensorFlow
\citep{abadi2016tensorflow, agrawal2019tensorflow}, JAX \citep{frostig18jax},
and Zygote \citep{innes2019zygote}. These AD tools are used widely, especially
to train machine learning models such as neural networks.

\paragraph{Optimization layers.} Implementing the derivative and its adjoint
of an optimization problem makes it possible to implement the problem as
a differentiable function in AD software. These differentiable solution
maps are sometimes called \emph{optimization layers} in the machine learning
community. Many specific optimization layers have been implemented, including
QP layers \citep{amos2017optnet}, convex optimization layers
\citep{agrawal2019differentiable}, and nonlinear program layers
\citep{gould2019deep}. Optimization layers have found several applications in,
\eg, computer graphics \citep{geng2020coercing}, control
\citep{agrawal2019cocp, de2018end, amos2018differentiable, barratt2019fitting},
data fitting and classification \citep{barratt2019least}, game playing
\citep{ling2018game}, and combinatorial tasks \citep{berthet2020learning}.

While some optimization layers are implemented by differentiating through each
step of an iterative algorithm (known as \emph{unrolling}), we emphasize that
in this paper, we differentiate through LLCPs \emph{analytically},
without unrolling an optimization algorithm, and without tracing each
step of the DSL.

\paragraph{Numerical solvers.} A numerical solver is an implementation of
an optimization algorithm, specialized to a specific subclass of optimization
problems. DSLs like CVXPY rewrite high-level descriptions of optimization
problems to the rigid low-level formats required by solvers. While 
some solvers have been implemented specifically for GPs
\citep[\S10.2]{boyd2007tutorial}, LLCPs (and GPs) can just as well be solved by
generic solvers for convex cone programs that support the exponential cone.
Our implementation reduces LLCPs to cone programs and solves them using SCS
\citep{odonoghue2016conic}, an ADMM-based solver for cone programs. In
principle, our method is compatible with other conic solvers as well, such as
ECOS \citep{domahidi2013ecos} and MOSEK \citep{mosek}.

\section{Disciplined geometric programming}\label{s-dgp}
The DGP ruleset for unparametrized LLCPs was given in \S\ref{sub-s-dsl}. Here,
we remark on the types of atoms under consideration, and
we extend DGP to include parameters as atoms. We then give several
examples of parametrized, DGP-compliant expressions.

\subsection{Atom library}
The class of LLCPs producible
using DGP depends on the atom library.
We make a few standard assumptions
on this library, limiting our attention to atoms that can be implemented in
a DSL for convex optimization. In particular, we assume that the log-log
transformation of each DGP atom (or its epigraph) can be represented in a
DCP-compliant fashion, using DCP atoms. For example, this means that if the
product is a DGP atom, then we require that the sum (which is its log-log
transformation) to be a DCP atom. In turn, we assume that the epigraph of each
DCP atom can be represented using the standard convex cones (\ie, the
zero cone, the nonnegative orthant, the second-order cone, the exponential
cone, and the semidefinite cone).

For simplicity, the reader may assume that the atoms under consideration are
the ones listed in the DGP tutorial at
\begin{center}
\texttt{https://www.cvxpy.org}.
\end{center} 
The subclass of LLCPs generated by these atoms is a superset of GGPs, because
it includes the product, sum, power, and maximum as atoms. It also includes
other basic functions, such as the ratio, difference, exponential, logarithm,
and entropy, and functions of elementwise positive matrices, such as the
spectral radius and resolvent.

\subsection{Parameters}\label{sub-s-parameters}
We extend the DGP ruleset to include parameters as atoms by defining the
curvature of parameters, and defining the curvature of a parametrized power
atom. Like unparametrized expressions, a parametrized expression is log-log
convex under DGP if can be generated by the composition rule. 

\paragraph{Curvature.} The curvature of a positive parameter is log-log affine.
Parameters that are not positive have unknown log-log curvature.

\paragraph{The power atom.} The power atom $f(x; a) = x^a$ is log-log
affine if the exponent $a$ is a fixed numerical constant, or if $a$ is 
parameter and the argument $x$ is not parametrized. If $a$
is a parameter, it need not be positive. The exponent $a$ is not an argument
of the power atom, \ie, DGP does not allow for the exponent to be a composition
of atoms. \\

These rules ensure that a parametrized DGP problem can be reduced to a
parametrized DCP problem, as explained in \S\ref{sub-s-canonicalization}. (This,
in turn, will simplify the calculation of the derivative of the solution map.)
These rules are similar to the rules for parameters from
\citep{agrawal2019differentiable}, in which a method for differentiating
through parametrized DCP problems was developed. They are are not too
restrictive; \eg, they permit taking the coefficients and exponents in a GP as
parameters. We now give several examples of the kinds of parametrized
expressions that can be constructed using DGP\@. 

\subsection{Examples}
\paragraph{Example 1.}
Consider a parametrized monomial
\[
c x_1^{a_1}x_2^{a_2}\cdots x_n^{a_n},
\]
where $x \in \reals^n_{++}$ is the variable and $c \in
\reals_{++}$, $a_1, \ldots, a_n \in \reals$ are parameters. This expression
is DGP-compliant. To see this, notice that each power expression $x_i^{a_i}$ is log-log affine,
since $a_i$ is a parameter and $x_i$ is not parametrized. Next, note that
the product of the powers is log-log affine, since the product of log-log affine
expressions is log-log affine. Finally, the parameter $c$ is log-log affine
because it is positive, so by the same reasoning the product of $c$ and
$x_1^{a_1}\cdots x_n^{a_n}$ is log-log affine as well.

On the other hand, if
$a_{n+1} \in \reals$ is an additional parameter, then
\[
{(c x_1^{a_1}x_2^{a_2}\cdots x_n^{a_n})}^{a_{n+1}},
\]
is not DGP-compliant, since $c x_1^{a_1}x_2^{a_2}\cdots
x_n^{a_n}$ and $a_{n+1}$ are both parametrized.

\paragraph{Example 2.}
Consider a parametrized posynomial, \ie, a sum of
monomials,
\[
\sum_{i=1}^{m} c_i x_1^{a_{i 1}}x_2^{a_{i2}}\cdots x_n^{a_{in}},
\]
where $x \in \reals^{n}_{++}$ is the variable, and $c \in \reals^{m}_{++}$,
$a_{ij} \in \reals$ ($i = 1, \ldots, m$, $j = 1, \ldots, n$) are parameters.
This expression is also DGP-compliant, since each term in the sum is a
parametrized monomial, which is log-log affine, and the sum of log-log affine
expressions is log-log convex.

\paragraph{Example 3.} The maximum of posynomials, parametrized
as in the previous examples, is log-log convex, since the maximum is a log-log
convex function that is increasing in each of its arguments.

\paragraph{Example 4.}
We can also give examples of parametrized expressions that do not involve monomials
or posynomials. In this example and the next one, a vector
or matrix expression is log-log convex (or log-log affine, or log-log concave)
if every entry is log-log convex (or log-log affine, or log-log concave).

\begin{itemize}
\item The expression $\exp(c \circ x)$, where $x \in \reals^n_{++}$ is the
variable and $c \in \reals^n_{++}$ is the parameter (and $\circ$ is the
elementwise product), is log-log convex, since $c \circ x$ is log-log
affine and $\exp$ is log-log convex; likewise, $\exp(c^Tx)$ is log-log convex,
since $c^Tx$ is log-log convex and $\exp$ is increasing.

\item The expression $\log(c \circ x)$, where $x \in \reals^n_{++}$ is the
variable and $c \in \reals^n_{++}$ is the parameter
is log-log concave, since $c \circ x$ is log-log affine and $\log$ is log-log
concave. However, $\log(c^Tx)$ does not have log-log curvature, since the
$\log$ atom is increasing and log-log concave but $c^Tx$ is log-log convex.
\end{itemize}

\paragraph{Example 5.} Finally, we give two examples involving functions of
matrices with positive entries.
\begin{itemize}
\item The spectral radius $\rho(X)$ of a matrix $X \in \reals^{n \times
n}_{++}$ with positive entries is a log-log convex function, increasing
in each entry of $X$. If $C \in \reals^{n \times n}_{++}$ is a parameter, then
$\rho(C \circ X)$ and $\rho(CX)$ are both log-log convex and DGP-compliant
($C \circ X$ is log-log affine, and $CX$ is log-log convex).

\item The atom $f(X) = {(I - X)}^{-1}$ is log-log convex (and increasing) in matrices
$X \in \reals^{n \times n}_{++}$ with $\rho(X) < 1$.
Therefore, if $X$ is a variable and $C \in \reals^{n \times n}_{++}$ is a
parameter, the expressions ${(I - C \circ X)}^{-1}$ and ${(I - C X)}^{-1}$ are
log-log convex and DGP-compliant.
\end{itemize}
We refer readers interested in these functions to \cite[\S2.4]{agrawal2019dgp}.

\section{The solution map and its derivative}\label{s-solution-map}
We consider a DGP-compliant LLCP, with variable $x \in \reals^{n}_{++}$
and parameter $\alpha \in \reals^{k}$. We assume throughout that
the solution map of the LLCP is single-valued, and we denote it by $\Smap
: \reals^{k} \to \reals^{n}_{++}$. There are several pathological
cases in which the solution map may not be differentiable; we simply
limit our attention to non-pathological cases, without explicitly
characterizing what those cases are. We leave a characterization of the
pathologies to future work. In this section we describe the form of the
implicit function $\Smap$, and we explain how to compute the derivative
operator $\DD \Smap$ and its adjoint $\DD^T \Smap$.

Because a DGP problem can be reduced to a DCP problem via the log-log
transformation, we can represent its solution map by the composition of
a map $\Cmap : \reals^{k} \to \reals^{p}$, which maps the parameters in the LLCP
to parameters in the DCP problem, the solution map $\phi : \reals^{p} \to
\reals^{m}$ of the DCP problem, and a map $\Rmap : \reals^{m} \to
\reals^{n}_{++}$ which recovers the solution to the LLCP from a solution to the
convex program. That is, we represent $\Smap$ as
\[
\Smap = \Rmap \circ \phi \circ \Cmap.
\]

In \S\ref{sub-s-canonicalization},  we describe the canonicalization map
$\Cmap$ and its derivative. When the conditions on parameters from
\S\ref{sub-s-parameters} are satisfied, the canonicalized parameters
$\Cmap(\alpha)$ (\ie, the parameters in the DCP program) satisfy the rules
introduced in \citep{agrawal2019differentiable}. This lets us use the method
from \citep{agrawal2019differentiable} to efficiently differentiate through the
log-log transformation of the LLCP, as we describe in \S\ref{sub-s-dcp}. In
\S\ref{sub-s-recovery}, we describe the recovery map $\Rmap$ and its
derivative.

Before proceeding, we make a few basic remarks on the
derivative.

\paragraph{The derivative.}
We calculate the derivative of $\Smap$ by calculating the derivatives of these
three functions, and applying the chain rule. Let $\beta = \Cmap(\alpha)$
and $\tilde{x}^{\star} = \phi(\beta)$. The derivative at $\alpha$ is just
\begin{equation*}\label{e-derivative}
\DD\Smap(\alpha) = \DD\Rmap(\tilde{x}^{\star}) \DD\phi(\beta) \DD\Cmap(\alpha),
\end{equation*}
and the adjoint of the derivative is
\begin{equation*}\label{e-adjoint}
\DD^T\Smap(\alpha) = \DD^T\Cmap(\alpha)  \DD^T\phi(\beta) \DD^T \Rmap(\tilde{x}^{\star}).
\end{equation*}
The derivative at $\alpha$, $\DD\Smap(\alpha)$, is a matrix in $\reals^{n \times k}$,
and its adjoint is its transpose.

\paragraph{Sensitivity analysis.} Suppose the parameter $\alpha$ is perturbed
by a vector $\dd\alpha \in \reals^k$ of small magnitude. Using the
derivative of the solution map, we can compute a first-order approximation of
the solution of the perturbed problem, \ie,
\[
\Smap(\alpha + \dd \alpha) \approx \Smap(\alpha) + \DD \Smap(\alpha)\dd\alpha.
\]
The quantity
\[
\DD \Smap(\alpha)\dd\alpha
\]
is an approximation of the change in the solution, due to the perturbation.

\paragraph{Gradient.}
Consider a function $f : \reals^n \to \reals$, and suppose we wish to compute
the gradient of the composition $f \circ \Smap$ at $\alpha$. By the chain rule, the
gradient is simply
\[
\nabla (f \circ \Smap)(\alpha) = \DD^T\Smap(\alpha) \dd x
\]
where $\dd x \in \reals^{n}$ is the gradient of $f$, evaluated at $\Smap(\alpha)$.
Notice that evaluating the adjoint of the derivative at a vector corresponds to
computing the gradient of a function of the solution. 
(In the machine learning community, this computation is known as backpropagation.)

\subsection{Canonicalization}\label{sub-s-canonicalization}
A DGP problem parametrized by $\alpha \in \reals^{k}$ can be
\emph{canonicalized}, or reduced, to an equivalent DCP problem parametrized by
$\beta \in \reals^{p}$. The canonicalization map $\Cmap : \reals^{k} \to
\reals^{p}$ relates the parameters $\alpha$ in the DGP problem to the
parameters $\beta$ in the DCP problem by $\Cmap(\alpha) = \beta$. In this
section, we describe the form of $\Cmap$, and explain why the problem produced by
canonicalization is DCP-compliant (with respect to the parametrized DCP
ruleset introduced in \cite{agrawal2019differentiable}).

The canonicalization of a parametrized DGP problem is the
same as the canonicalization of an unparametrized DGP problem in which the 
parameters have been replaced by constants. A DGP expression can be thought
of an expression tree, in which the leaves are variables, constants,
or parameters, and the root and inner nodes are atomic functions. The children
of a node are its arguments. Canonicalization recursively replaces each 
expression with its log-log transformation, or the log-log transformation of
its epigraph \citep[\S4.1]{agrawal2019dgp}. For example, positive variables are
replaced with unconstrained variables and products are replaced with sums.

Parameters appearing as arguments to an atom are replaced with
their logs. Parameters appearing as exponents in power atoms, however,
enter the DCP problem unchanged; the expression $x^a$ is canonicalized to
$a F(u)$, where $F(u)$ is the log-log transformation of the expression $x$.
In particular, for $\beta = \Cmap(\alpha)$,
for each $i = 1, \ldots, p$, there exists $j \in \{1, \ldots, k\}$ such that
either $\beta_i = \log(\alpha_j)$ or $\beta_i = \alpha_j$. The derivative
$\mathsf{D}\Cmap(\alpha) \in \reals^{p \times k}$ is therefore
easy to compute. Its entries are given by
\[
{\mathsf{D}\Cmap(\alpha)}_{ij} =
\begin{cases}
1 & \beta_i = \alpha_j \\
1/\alpha_j & \beta_i = \log(\alpha_j) \\
0 & \text{otherwise.}
\end{cases}
\]
for $i = 1, \ldots, p$ and $j = 1, \ldots, k$.

\paragraph{DCP compliance.} When the DGP problem is not parametrized, the
convex optimization problem emitted by canonicalization is DCP-compliant
\citep{agrawal2019dgp}. When the DGP problem is parametrized, it turns out
that the emitted convex optimization problem satisfies the DCP ruleset for
parametrized problems, given in \citep[\S4.1]{agrawal2019differentiable}. In
DCP, parameters are affine (just as parameters are log-log affine in DGP);
additionally, the product $xy$ is affine if either $x$ or $y$ is a numerical
constant, $x$ is a parameter and $y$ is not parametrized, or $y$ is a parameter
and $x$ is not parametrized. The restriction on the power atom in DGP
ensures that all products appearing in the emitted convex optimization
problem are affine under DCP\@. DCP-compliance of the remaining expressions
follows from the assumptions on the atom library, which guarantee that the
log-log transformation of a DGP expression is DCP-compliant. Therefore, the
parametrized convex optimization problem is DCP-compliant. (In
the terminology of \citep{agrawal2019differentiable}, the problem is a
\emph{disciplined parametrized program}.)

\subsection{The convex optimization problem}\label{sub-s-dcp}
The convex optimization problem emitted by canonicalization has the variable
$\tilde x \in \reals^{m}$, with $m \geq n$.
The solution map $\phi$ of the DCP problem maps the
parameter $\beta \in \reals^{p}$ to the solution $\tilde{x}^\star \in
\reals^{m}$. Because the problem is DCP-compliant, to
compute its derivative $\DD\phi(\beta)$, we can simply use the method from
\citep{agrawal2019differentiable}.

In particular, $\phi$ can be represented in \emph{affine-solver-affine}
form: it is the composition of an affine map from parameters in the DCP
problem to the problem data of a convex cone program; the solution map of a
convex cone program; and an affine map from the solution of the convex cone
program to the solution of the DCP problem. The affine maps and their
derivatives can be evaluated efficiently, since they can be
represented as sparse matrices \citep[\S4.2, 4.4]{agrawal2019differentiable}.
The cone program can be solved using standard algorithms for conic
optimization, and its derivative can be computed using the
method from \citep{diffcp2019}; the latter involves computing certain
projections onto cones, their derivatives, and solving a least-squares problem.

\subsection{Solution recovery}\label{sub-s-recovery}
Let $\tilde{x} \in \reals^m$ be the variable in the DCP problem, and
partition $\tilde{x}$ as $(\hat x, s) \in \reals^{n \times (m-n)}$. Here,
$\hat x$ is the elementwise log of the variable $x$ in the DGP problem
and $s$ is a slack variable involved in graph implementations of DGP
atoms. If $(\hat{x}^{\star}, s^{\star})$ is optimal for the DCP problem,
then $\exp(\hat{x}^{\star})$ is optimal for the DGP problem (the
exponentiation is meant elementwise) \citep[\S4.2]{agrawal2019dgp}. Therefore,
the recovery map $\Rmap : \reals^{m} \to \reals^{n}$ is given by
\[
\Rmap(\tilde x) = \exp(\hat x).
\]
The entries of its derivative are simply given by
\[
\DD\Rmap(\tilde x)_{ij} =
\begin{cases}
\exp(\hat{x}_i) & i = j \\
0  & \textnormal{otherwise,}
\end{cases}
\]
for $i = 1, \ldots, n$, $j = 1, \ldots, m$.

\section{Implementation}\label{s-implementation}
We have implemented the derivative and adjoint derivative of LLCPs
as abstract linear operators in CVXPY, a Python-embedded modeling
language for convex optimization and log-log convex optimization
\citep{diamond2016cvxpy, agrawal2018rewriting, agrawal2019dgp}. With our
software, users can differentiate through any parametrized LLCP produced via
the DGP ruleset, using the atoms listed at
\begin{center}
\texttt{https://www.cvxpy.org}.
\end{center}
Additionally, we provide differentiable PyTorch and TensorFlow layers for
LLCPs in CVXPY Layers, available at
\begin{center}
\texttt{https://www.github.com/cvxgrp/cvxpylayers}.
\end{center}

We now remark on a few aspects of our implementation, before presenting usage
examples in \S\ref{sub-s-example}.

\paragraph{Caching.}
In our implementation, the first time a parametrized DGP problem is solved, we
compute and cache the canonicalization map $\Cmap$, the parametrized
DCP problem, and the recovery map $\Rmap$. On subsequent solves, instead of
re-canonicalizing the DGP problem, we simply evaluate $\Cmap$ at the parameter
values and update the parameters in the DCP problem in-place. The
affine maps involved in the canonicalization and solution recovery of the DCP
problem are also cached after the first solve, as in
\citep{agrawal2019differentiable}. This means after an initial ``compilation'',
the overhead of the DSL is negligible compared to the time spent in the
numerical solver.

\paragraph{Derivative computation.} CVXPY (and CVXPY Layers) represent the
derivative and its adjoint abstractly, letting users evaluate them at vectors.
If $\alpha \in \reals^{k}$ is the parameter, and $\dd\alpha \in \reals^{k}$,
$\dd x \in \reals^{n}$ are perturbations, users may compute
$\DD\Smap(\alpha)(\dd\alpha)$ and $\DD^T\Smap(\alpha)(\dd x)$. In particular,
we do not materialize the derivative matrices. Because the action of
the DSL is cached after the first compilation, we can compute these operations
efficiently, without tracing each instruction executed by the DSL\@.

\subsection{Hello world}\label{sub-s-example}
Here, we present a basic example of how
to use CVXPY to specify a parametrized and solve DGP problem, and how
to evaluate its derivative and adjoint. This example is only meant to illustrate
the usage of our software; a more interesting example is presented in
\S\ref{s-example}.

Consider the following code:
\begin{lstlisting}
import cvxpy as cp

x = cp.Variable(pos=True)
y = cp.Variable(pos=True)
z = cp.Variable(pos=True) 

a = cp.Parameter(pos=True)
b = cp.Parameter(pos=True)
c = cp.Parameter()

objective_fn = 1/(x*y*z)
objective = cp.Minimize(objective_fn)
constraints = [a*(x*y + x*z + y*z) <= b, x >= y**c]
problem = cp.Problem(objective, constraints)

print(problem.is_dgp(dpp=True))
\end{lstlisting} 
This code block constructs an LLCP \texttt{problem}, with
three scalar variables, $x, y, z \in \reals_{+}$. Notice that the variables
are declared as positive, with \texttt{pos=True}. The objective is to minimize
the reciprocal of the product of the variables, which is log-log affine.
There are three parameters, \texttt{a}, \texttt{b}, and \texttt{c}, two of which
are declared as positive.
The variables are constrained so that \texttt{x} is at least \texttt{y**c}
(\ie, \texttt{y} raised to the power \texttt{c}); notice that this constraint
is DGP-compliant, since the power atom is log-log affine and \texttt{x} is
log-log affine. Additionally, the posynomial \texttt{a*(x*y + x*z + y*z)}
is constrained to be no larger than the parameter \texttt{b}; the posynomial
is log-log convex and DGP-compliant, since the parameters are positive.
The penultimate line constructs the problem, and the last line checks whether
the problem is DGP; the keyword argument \texttt{dpp=True} tells CVXPY that it
should use the rules involving parameters introduced in
\S\ref{sub-s-parameters}. As expected, the output of this program is the string
\texttt{True}.

\paragraph{Solving the problem.} The problem constructed above can be solved
in one line, after setting the values of the parameters, as below.
\begin{lstlisting}
a.value = 2.0
b.value = 1.0
c.value = 0.5
problem.solve(gp=True, requires_grad=True)
\end{lstlisting}
The keyword argument \texttt{gp=True} tells CVXPY to parse the problem
using DGP, and the keyword argument \texttt{requires\_grad=True} will let
us subsequently evaluate the derivative and its adjoint. After calling \texttt{problem.solve},
the optimal values of the variables are stored in the \texttt{value} attribute,
that is,
\begin{lstlisting}
print(x.value)
print(y.value)
print(z.value)
\end{lstlisting}
prints 
\begin{lstlisting}
0.5612147353889386
0.31496200373359456
0.36892055859991446
\end{lstlisting}
(and the optimal value of the problem is stored in \texttt{problem.value}).

\paragraph{Sensitivity analysis.}
Suppose we perturb the parameter vector $\alpha$ by
a vector $\dd \alpha$ of small magnitude. We can approximate the change $\Delta$ in the
solution due to the perturbation using the derivative of the solution map,
as
\[
\Delta = \Smap(\alpha + \dd \alpha) - \Smap(\alpha) \approx \DD\Smap(\alpha)\dd\alpha.
\]

We can compute this quantity in CVXPY\@. For our running example, partition
the perturbation as
\[
\dd\alpha = \begin{bmatrix}
\dd a \\
\dd b \\
\dd c
\end{bmatrix}.
\]
To approximate the change in the optimal values for the variables $x$, $y$,
and $z$, we set the \texttt{delta} attributes on the parameters and then
call the \texttt{derivative} method.
\begin{lstlisting}
a.delta = da
b.delta = db
c.delta = dc
problem.derivative()
\end{lstlisting}
The \texttt{derivative} method populates the \texttt{delta} attributes of the
variables in the problem as a side-effect. Say we set \texttt{da}, \texttt{db},
and \texttt{dc} to 1e-2. Let $\hat x$, $\hat y$, and $\hat z$ be the
first-order approximations of the solution to the perturbed problem; we
can compare these to the actual solution, as follows.
\begin{lstlisting}
x_hat = x.value + x.delta
y_hat = y.value + y.delta
z_hat = z.value + z.delta

a.value += da
b.value += db
c.value += dc
problem.solve(gp=True)

print('x: predicted {0:.5f} actual {1:.5f}'.format(x_hat, x.value))
print('y: predicted {0:.5f} actual {1:.5f}'.format(y_hat, y.value))
print('z: predicted {0:.5f} actual {1:.5f}'.format(z_hat, z.value))
\end{lstlisting}

\begin{lstlisting}
x: predicted 0.55729 actual 0.55732
y: predicted 0.31783 actual 0.31781
z: predicted 0.37179 actual 0.37178
\end{lstlisting}

\paragraph{Gradient.}
We can compute the gradient of a function of the solution with respect to
the parameters, using the adjoint of the derivative of the solution map. Let
$\alpha = (a, b, c)$ be the parameters in our problem, and let $x(\alpha)$,
$y(\alpha)$, and $z(\alpha)$ denote the optimal variable values for our
problem, so that
\[
\Smap(\alpha) = \begin{bmatrix}x(\alpha) \\ y(\alpha) \\ z(\alpha)\end{bmatrix},
\]
where $\Smap$ is the solution map of our optimization problem. Let $f :
\reals^{3} \to \reals$, and suppose we wish to compute the gradient of the
composition $f \circ \Smap$ at $\alpha$. By the chain rule,
\[
\nabla f(\Smap(\alpha)) = \DD^T\Smap(\alpha) \begin{bmatrix}
\dd x \\
\dd y \\
\dd z, 
\end{bmatrix}
\]
where $\dd x, \dd y, \dd z$ are the partial derivatives of $f$ with respect to
its arguments.

We can compute the gradient in CVXPY\@. Below, \texttt{dx}, \texttt{dy}, and
\texttt{dz} are numerical constants, corresponding
to $\dd x, \dd y$, and $\dd z$.
\begin{lstlisting}
x.gradient = dx
y.gradient = dy
z.gradient = dz
problem.backward()
\end{lstlisting}
The \texttt{backward} method populates the \texttt{gradient} attributes on the
parameters. If left uninitialized, the \texttt{gradient} attributes on the
variables
default to $1$, corresponding to taking $f$ to be the sum function. The
gradient of a scalar-valued function $f$ with respect to the solution (the
values $\dd x, \dd y, \dd z$) may be computed manually, or using software
for automatic differentiation. 

As an example, suppose $f$ is the function
\[
\begin{bmatrix} x \\ y \\ z \end{bmatrix} \mapsto \frac{1}{2}(x^2 + y^2 + z^2),
\]
so that $\dd x = x$, $\dd y = y$, and $\dd z = z$.
Let $\dd\alpha = \nabla f(\Smap(\alpha))$, and say we subtract $\eta \dd\alpha$
from the parameter, where $\eta$ is a small positive number, such as $0.5$.
Using the following code, we can compare $f(\Smap(\alpha - \dd\alpha))$ with the
value predicted by the gradient, \ie
\[
f(\Smap(\alpha - \eta\dd\alpha)) \approx f(\Smap(\alpha)) - \eta\dd\alpha^T \dd\alpha.
\]
\begin{lstlisting}
def f(x, y, z):
    return 1/2*(x**2 + y**2 + z**2)

original = f(x, y, z).value

x.gradient = x.value
y.gradient = y.value
z.gradient = z.value
problem.backward()

eta = 0.5
dalpha = cp.vstack([a.gradient, b.gradient, c.gradient])
predicted = float((original - eta*dalpha.T @ dalpha).value)

a.value -= eta*a.gradient
b.value -= eta*b.gradient
c.value -= eta*c.gradient
problem.solve(gp=True)
actual = f(x, y, z).value

print('original {0:.5f} predicted {1:.5f} actual {2:.5f}'.format(
       original, predicted, actual))
\end{lstlisting}

\begin{lstlisting}
original 0.27513 predicted 0.22709 actual 0.22942
\end{lstlisting}
\paragraph{CVXPY Layers.}
We have implemented support for solving and differentiating through LLCPs
specified with CVXPY in CVXPY Layers, which provides PyTorch and TensorFlow
wrappers for our software. This makes it easy to use automatic differentiation
to compute the gradient of a function of the solution. For example, the above
gradient calculation can be done in PyTorch, with the following code.
\begin{minipage}{\linewidth}
\begin{lstlisting}
from cvxpylayers.torch import CvxpyLayer
import torch

layer = CvxpyLayer(problem, parameters=[a, b, c],
                   variables=[x, y, z], gp=True)
a_tch = torch.tensor(2.0, requires_grad=True)
b_tch = torch.tensor(1.0, requires_grad=True)
c_tch = torch.tensor(0.5, requires_grad=True)

x_star, y_star, z_star = layer(a_tch, b_tch_c_tch)
sum_of_solution = x_star + y_star + z_star
sum_of_solution.backward()
\end{lstlisting}
\end{minipage}
The PyTorch method \texttt{backward} computes the gradient of the sum of the
solution, and populates the \texttt{grad} attribute on the PyTorch tensors which
were declared with \texttt{requires\_grad=True}. Of course, we could just as
well replace the sum operation on the solution with another scalar-valued
operation.

\subsection{Performance}
Here, we report the time it takes our software to parse, solve, and differentiate through
a DGP problem of modest size. The problem under consideration is a GP,
\begin{equation}
\begin{array}{ll}
\mbox{minimize} & \prod_{j=1}^{n}x_i^{A_{1j}} \\
\mbox{subject to} & \sum_{i=1}^{m} c_j  \prod_{j=1}^{n}x_i^{A_{ij}}  \leq 1 \\
& l \leq x \leq u,
\end{array}\label{e-benchmark}
\end{equation}
with variable $x \in \reals^{n}_{++}$ and parameters $A \in \reals^{m \times
n}$, $c \in \reals^{m}_{++}$, $l \in \reals^{n}_{++}$, and $u \in
\reals^{n}_{++}$.

We solve a specific numerical instance, with $n = 5000$ and $m=3$,
corresponding a problem with $5000$ variables and $25003$ parameters. We solve
the problem using SCS \citep{odonoghue2016conic, odonoghue2017scs}, and use
\texttt{diffcp} to compute the derivatives \citep{diffcp2019, diffcpsoftware}.
In table~\ref{table-timings}, we report the mean $\mu$ and standard deviation
$\sigma$ of the wall-clock times for the \texttt{solve}, \texttt{derivative},
and \texttt{backward} methods, over 10 runs (after performing a warm-up
iteration). These experiments were conducted on a standard laptop, with 16 GB
of RAM and a 2.7 GHz Intel Core i7 processor. We break down the report into the
time spent in CVXPY and the time spent in the numerical solver; the total
wall-clock time is the sum of these two quantities.

\begin{table}[]
\centering
\begin{tabular}{lllll}
\hline
& $\mu_{\mathrm{CVXPY}}$ & $\sigma_{\mathrm{CVXPY}}$ & $\mu_{\mathrm{solver}}$ & $\sigma_{\mathrm{solver}}$ \\ \hline \hline \\
\texttt{solve} & 12.54 & 1.03 & 2753.32 & 87.92 \\
\texttt{derivative} & 6.27 & 0.83 & 2.94  & 0.29  \\
\texttt{backward} & 37.2 & 0.88 & 19.89 & 0.59
\end{tabular}
\caption{Timings for problem~\ref{e-benchmark}, in ms.}\label{table-timings}
\end{table}

Because our implementation caches compact representations of the
canonicalization map and recovery map, instead of tracing the entire execution
of the DSL, the overhead of CVXPY is negligible compared to the time spent in
the numerical solver. For the \texttt{solve} method, the time spent in CVXPY
--- which maps the parameters in the LLCP to the parameters in a convex cone
program, and retrieves a solution of the LLCP from a solution of the cone
program --- is roughly two orders of magnitude less than the time spent in the
numerical solver, accounting for just 0.45\% of the total wall-clock time.
Similarly, computing the derivative (and its adjoint) is also about two orders of
magnitude faster than solving the problem. 

\section{Examples}\label{s-example}
In this section, we present two illustrative examples. The first example uses
the derivative of an LLCP to analyze the design of a queuing system. The second
example uses the adjoint of the derivative, training an LLCP as an optimization
layer for a synthetic structured regression task. 

The code for our examples are available online, at
\begin{center}
\texttt{https://www.cvxpy.org/examples/index.html}.
\end{center}

\subsection{Queuing system}
We consider the optimization of a (Markovian) \emph{queuing system}, with $N$
queues. A queuing system is a collection of queues, in which queued items
wait to be served; the queued items might be threads in an operating system,
or packets in an input or output buffer of a networking system. A natural goal
to minimize the service load of the system, given constraints on various
properties of the queuing system, such as limits on the maximum delay or latency.
In this example, we formulate this design problem as an LLCP, and compute the
sensitivity of the design variables with respect to the parameters.
The queuing system under consideration here is known as an $M/M/N$ queue, in
Kendall's notation \citep{kendall1953stochastic}. Our formulation follows
\citep{chiang2002efficient}.

We assume that items arriving at the $i$th queue are generated by a
Poisson process with rate $\lambda_i$, and that the service times for the $i$th
queue follow an exponential distribution with parameter $\mu_i$, for $i=1,
\ldots, N$. The \emph{service load} of the queuing system is a
function $\ell : \reals^{N}_{++} \times \reals^{N}_{++} \to \reals^{N}_{++}$
of the arrival rate vector $\lambda$ and the service rate vector
$\mu$, with components
\[
\ell_i(\lambda, \mu) = \frac{\mu_i}{\lambda_i}, \quad i=1, \ldots, N.
\]
(This is the reciprocal of the traffic load, which is usually denoted by $\rho$.)
Similarly, the queue occupancy, the average delay, and the total delay
of the system are (respectively) functions $q$, $w$, and $d$ of $\lambda$
and $\mu$, with components
\[
q_i(\lambda, \mu) =
\frac{\ell_i(\lambda, \mu)^{-2}}{1 - \ell_i(\lambda, \mu)^{-1}}, \quad
w_i(\lambda, \mu) = \frac{q_i(\lambda, \mu)}{\lambda_i} + \frac{1}{\mu_i}, \quad
d_i(\lambda, \mu) = \frac{1}{\mu_i - \lambda_i}
\]
These functions have domain $\{(\lambda, \mu) \in \reals^{N}_{++} \times
\reals^{N}_{++} \mid \lambda < \mu \}$, where the inequality
is meant elementwise. The queuing system has limits on the queue occupancy,
average queuing delay, and total delay, which must satisfy
\[
q(\lambda, \mu) \leq q_{\max}, \quad w(\lambda, \mu) \leq w_{\max}, \quad d(\lambda, \mu) \leq d_{\max},
\]
where $q_{\max}$, $w_{\max}$, and $d_{\max} \in \reals^{N}_{++}$ are
parameters and the inequalities are meant elementwise. Additionally, the arrival
rate vector $\lambda$ must be at least $\lambda_{\mathrm{min}} \in \reals^{N}_{++}$, and the sum of
the service rates must be no greater than $\mu_{\max} \in \reals_{++}$. 

Our design problem is to choose the arrival rates and service times to
minimize a weighted sum of the service loads, $\gamma^T \ell(\lambda, \mu)$,
where $\gamma \in \reals^{N}_{++}$ is the weight vector, while satisfying the
constraints. The problem is
\[
\begin{array}{ll}
\mbox{minimize} & \gamma^T \ell(\lambda, \mu) \\
\mbox{subject to}
& q(\lambda, \mu) \leq q_{\max} \\
& w(\lambda, \mu) \leq w_{\max} \\
& d(\lambda, \mu) \leq d_{\max} \\
& \lambda \geq \lambda_{\mathrm{min}}, \quad
\sum_{i=1}^{N} \mu_i \leq \mu_{\max}.
\end{array}\
\]
Here, $\lambda, \mu \in \reals^{N}_{++}$ are the
variables and $\gamma, q_{\max}, w_{\max}, d_{\max},
\lambda_{\mathrm{min}} \in \reals^{N}_{++}$ and $\mu_{\max} \in \reals_{++}$
are the parameters. This problem is an LLCP\@.
The objective function is a posynomial, as is the constraint function $w$. The functions
$d$ and $q$ are not posynomials, but they are log-log convex; log-log convexity
of $d$ follows from the composition rule, since the
function $(x, y) \mapsto y - x$ is log-log concave (for $0 < x < y$), and the ratio $(x, y) \mapsto
x/y$ is log-log affine and decreasing in $y$. By a similar argument,
$q$ is also log-log convex.

\paragraph{Numerical example.}
We specify a specific numerical instance of this problem in CVXPY using DGP,
with $N=2$ queues and parameter values
\[
\gamma = \begin{bmatrix} 1 \\ 2\end{bmatrix}, ~
q_{\max} = \begin{bmatrix}4 \\ 5\end{bmatrix}, ~
w_{\max} = \begin{bmatrix}2.5 \\ 3\end{bmatrix}, ~
d_{\max} = \begin{bmatrix}2 \\ 2\end{bmatrix}, ~
\lambda_{\mathrm{min}} = \begin{bmatrix}0.5\\ 0.8\end{bmatrix}, ~
\mu_{\max} = 3.
\]

We first solve the problem using SCS \citep{odonoghue2016conic, odonoghue2017scs}, obtaining the optimal values
\[
\lambda^{\star} = \begin{bmatrix}0.828 \\ 1.172\end{bmatrix}, \quad
\mu^{\star} = \begin{bmatrix}1.328 \\ 1.672\end{bmatrix}.
\]
Next, we perform a basic sensitivity analysis
by perturbing the parameters by one percent of their values, and computing
the percent change in the optimal variable values predicted by a first-order
approximation; we use CVXPY and the \texttt{diffcp} package \citep{diffcp2019}
to differentiate through the LLCP\@. We then compare the predicted change with
the true change by re-solving the problem at the perturbed values. The
predicted and true changes are
\[
\Delta \lambda_{\mathrm{pred}} = \begin{bmatrix}+2.3\% \\ +1.8\%\end{bmatrix}, ~
\Delta \lambda_{\mathrm{true}} = \begin{bmatrix}+2.0\% \\ +2.0\%\end{bmatrix}, ~
\Delta \mu_{\mathrm{pred}} = \begin{bmatrix}+1.1\% \\ +0.9\%\end{bmatrix}, ~
\Delta \mu_{\mathrm{true}} = \begin{bmatrix}+0.9\% \\ +1.1\%\end{bmatrix},
\]

\begin{table}
\centering
\begin{tabular}{llll}
\hline
& $d_{\max}$ & $\mu_{\max}$ & $\gamma$ \\ \hline \hline \\
$\lambda^\star$ &
    $\begin{bmatrix} -0.028 & 0.30 \\ 0.028 & -0.052\end{bmatrix}$ &
    $\begin{bmatrix} 0.46 \\ 0.54\end{bmatrix}$ &
    $\begin{bmatrix} 0.34 & -0.17 \\ -0.34 & 0.17\end{bmatrix}$ \\ \\
$\mu^\star$ &
    $\begin{bmatrix} -0.28 & 0.30 \\ 0.28 & -0.30\end{bmatrix}$ &
    $\begin{bmatrix} 0.46 \\ 0.54\end{bmatrix}$ &
    $\begin{bmatrix} 0.34 & -0.17 \\ -0.34 & 0.17\end{bmatrix}$ \\ \\
$\ell(\lambda^\star, \mu^\star)$ &
    $\begin{bmatrix} -0.28 & -0.22 \\ -0.10 & -0.20\end{bmatrix}$ &
    $\begin{bmatrix} -0.33 \\ -0.20\end{bmatrix}$ &
    $\begin{bmatrix} -0.24 & 0.12 \\ 0.12 & -0.061\end{bmatrix}$ \\
\end{tabular}
\caption{Derivatives of $\lambda^\star$, $\mu^\star$, and $\ell$ with respect to
$d_{\max}$, $\mu_{\max}$, and $\gamma$.}\label{table-d-var}
\end{table}
To examine the sensitivity of the solution to the individual
parameters, we compute the derivative of the variables with respect to
the parameters. The derivatives with respect to
$w_{\max}$, $q_{\max}$, and $\lambda_{\mathrm{min}}$ are
essentially $0$ (on the order of 1e-10), meaning that these parameters
can be changed slightly without affecting the solution. The derivatives
of the solution with respect to $d_{\max}$, $\mu_{\max}$, and $\gamma$
are given in table~\ref{table-d-var}.
While the solution is insensitive to small changes to the
limits on the queue occupancy, average queuing delay, and arrival rate,
it is highly sensitive to the limits on the total delay and
service rate, and to the weighting vector $\gamma$.

Finally, table~\ref{table-d-var} also lists the derivative of the service load
$\ell(\lambda^\star, \mu^\star)$ with respect to the parameters $d_{\max},
\mu_{\max}$, and $\gamma$. The table suggests that increasing the limits on the
total delay $d_{\max}$ and service rate $\mu_{\max}$ would decrease the service
loads on both queues, especially the first queue.

\subsection{Structured prediction}
In this example, we fit a regression model to structured data, using an LLCP\@.
The training dataset $\mathcal D$ contains $N$ input-output pairs $(x, y)$,
where $x \in \reals^{n}_{++}$ is an input and $y \in \reals^{m}_{++}$ is an
outputs. The entries of each output $y$ are sorted in ascending order, meaning
$y_1 \leq y_2 \leq \cdots y_m$.

Our regression model $\phi : \reals^{n}_{++} \to \reals^{m}_{++}$ takes as
input a vector $x \in \reals^{n}_{++}$, and solves an LLCP to produce a
prediction $\hat y \in \reals^{m}_{++}$. In particular, the solution of the
LLCP is the model's prediction. The model is of the form
\begin{equation}
\begin{array}{lll}
\phi(x) = &
\mbox{argmin} & \ones^T (z/y + y / z) \\
& \mbox{subject to} &  y_i \leq y_{i+1}, \quad i=1, \ldots, m-1 \\
&& z_i = c_i x_1^{A_{i1}}x_2^{A_{i2}}\cdots x_n^{A_{in}}, \quad i = 1, \ldots, m.
\end{array}\label{e-model}
\end{equation}
Here, the minimization is over $y \in \reals^{m}_{++}$ and an auxiliary
variable $z \in \reals^{m}_{++}$, $\phi(x)$ is the optimal value of $y$, and
the parameters are $c \in \reals^{m}_{++}$ and $A \in \reals^{m \times n}$. The
ratios in the objective are meant elementwise, as is the inequality $y \leq z$, and
$\ones$ denotes the vector of all ones. Given a vector $x$, this model finds a
sorted vector $\hat y$ whose entries are close to monomial functions of $x$
(which are the entries of $z$), as measured by the fractional error.

The training loss
$\mathcal{L}(\phi)$ of the model on the training set is the mean squared loss
\[
\mathcal{L}(\phi) = \frac{1}{N}\sum_{(x, y) \in \mathcal D} \norm{y - \phi(x)}_2^2.
\]
We emphasize that $\mathcal{L}(\phi)$ depends on $c$ and $A$.
In this example, we fit the parameters $c$ and $A$ in the LLCP~(\ref{e-model})
to minimize the training loss $\mathcal{L}(\phi)$.

\paragraph{Fitting.} We fit the parameters by an iterative projected
gradient descent method on $\mathcal L(\phi)$. In each iteration, we first
compute predictions $\phi(x)$ for each input in the training set; this requires
solving $N$ LLCPs. Next, we evaluate the training loss $\mathcal L(\phi)$. To
update the parameters, we compute the gradient $\nabla \mathcal L(\phi)$ of the
training loss with respect to the parameters $c$ and $A$. This
requires differentiating through the solution map of the LLCP~(\ref{e-model}).
We can compute this gradient efficiently, using the adjoint of the solution
map's derivative, as described in \S\ref{s-solution-map}. Finally, we subtract
a small multiple of the gradient from the parameters. Care must be taken to
ensure that $c$ is strictly positive; this can be done by clamping the entries
of $c$ at some small threshold slightly above zero. We run this method for
a fixed number of iterations.

\begin{figure}
\includegraphics{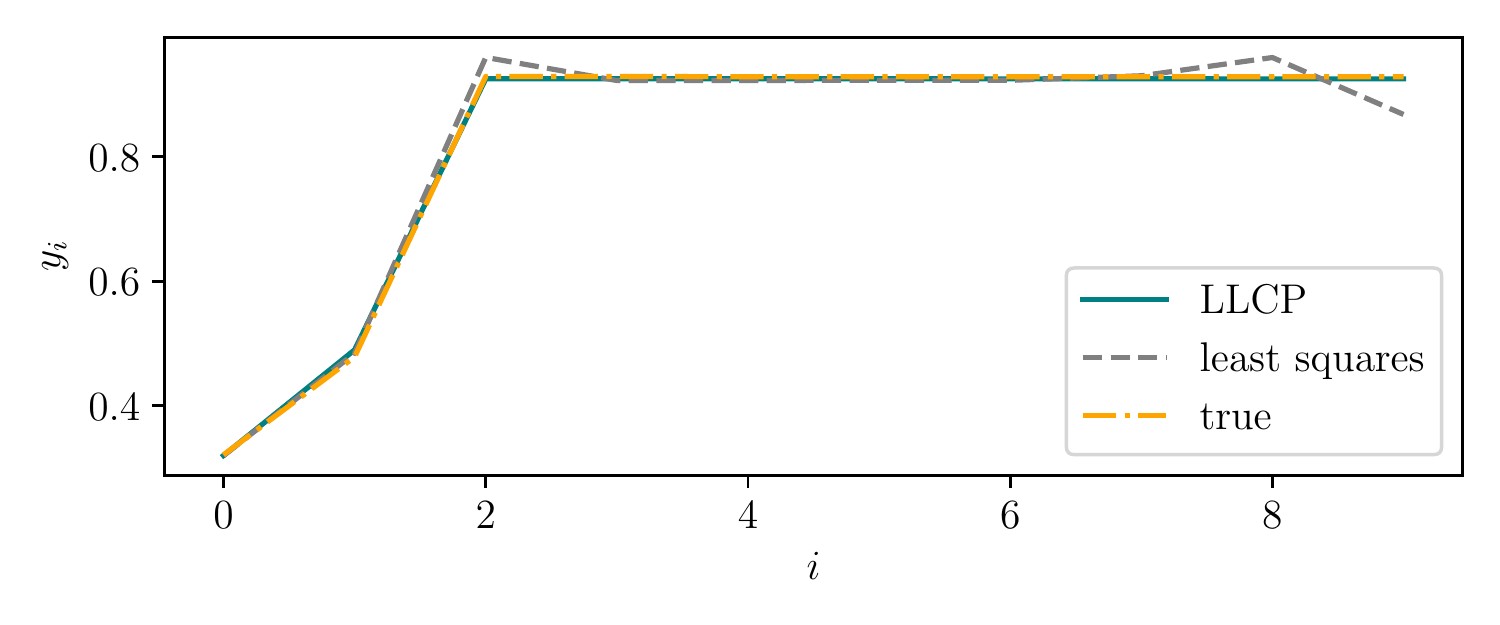}
\includegraphics{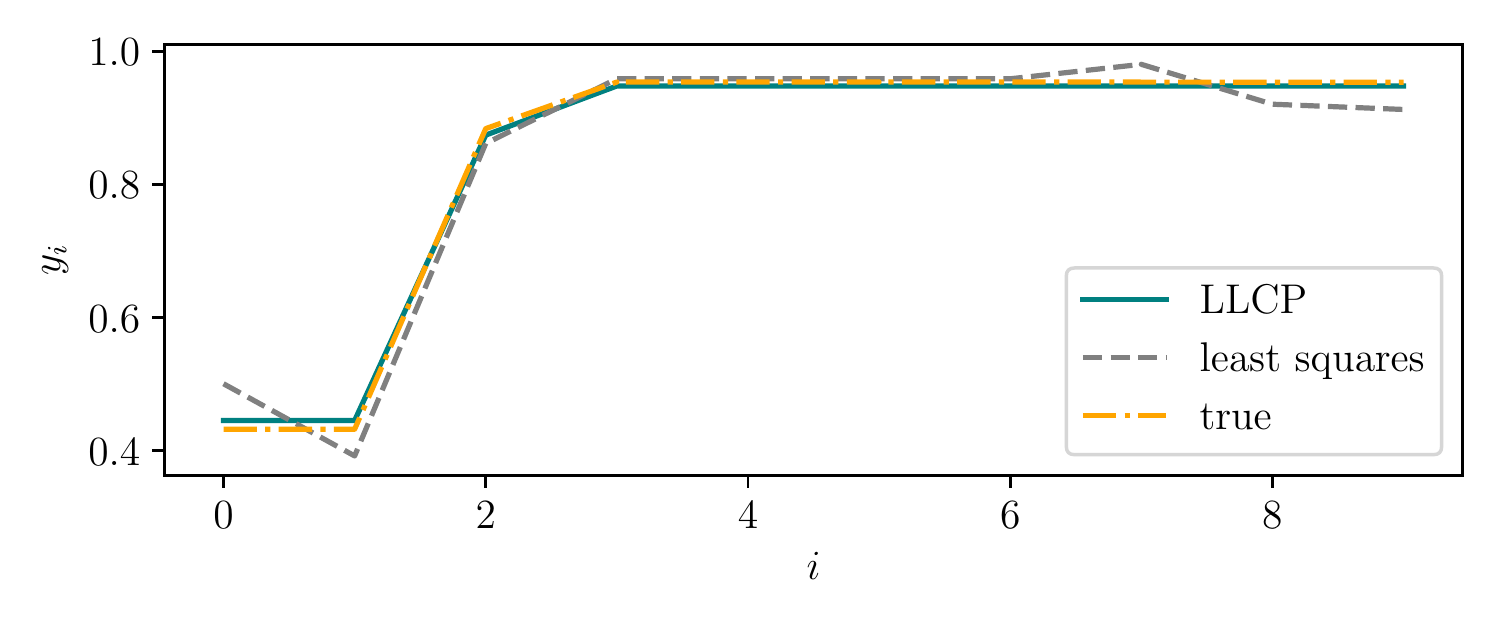}
\caption{Sample predictions and true output.}\label{f-sample}
\end{figure}

\paragraph{Numerical example.} We consider a specific numerical example, with
$N = 100$ training pairs, $n = 20$, and $m = 10$. The inputs were chosen according to
\[
x = \exp(\tilde x), \quad \tilde x \sim \mathcal{N}(0, I).
\]
We generated true parameter values
$A^\star \in \reals^{m \times n}$ (with entries sampled from a normal
distribution with zero mean and standard deviation $0.1$) and
$c^\star \in \reals^{m}_{++}$ (with entries set to the absolute value of
samples from a standard normal). The outputs were generated by
\[
y = \phi(x + \exp(v); A^{\star}, c^{\star}), \quad v \sim \mathcal{N}(0, I).
\]
We generated a held-out validation set of $50$ pairs, using the same
true parameters.

We implemented the LLCP~(\ref{e-model}) in CVXPY, and used PyTorch and
CVXPY Layers to train it using our gradient method. To initialize
the parameters $A$ and $c$, we computed a least-squares monomial fit to
the training data to obtain $A^{\mathrm{lstsq}}$ and $c^{\mathrm{lstsq}}$,
via the method described in \citep[\S8.3]{boyd2007tutorial}. From this
initialization, we ran 10 iterations of our gradient method. Each iteration,
which requires solving and differentiating through 150 LLCPs (100 for the
training data, and 50 for logging the validation error), took roughly 10
seconds on a 2012 MacBook Pro with 16 GB of RAM and a 2.7 GHz Intel Core i7
processor. 

The least-squares fit has a validation error of 0.014. The LLCP, with
parameters $A = A^{\mathrm{lstsq}}$ and $c = c^{\mathrm{lstsq}}$, has a
validation error of 0.0081, which is reduced to to 0.0077 after training.
Figure~\ref{f-sample} plots sample predictions of the LLCP and the
least-squares fit on a validation input, as well as the true output. The LLCP's
prediction is monotonic, while the least squares prediction is not. 

\section{Acknowledgments}
We thank Shane Barratt, who suggested using an optimization layer
to regress on sorted vectors, and Steven Diamond and Guillermo Angeris, for
helpful discussions related to data fitting and the derivative implementation.
Akshay Agrawal is supported by a Stanford Graduate Fellowship.

\clearpage
\printbibliography{}
\end{document}

%% file: defs.tex
\newcommand{\eg}{{\it e.g.}}
\newcommand{\ie}{{\it i.e.}}

\newcommand{\ones}{\mathbf 1}

\newcommand{\reals}{{\mbox{\bf R}}}

\newcounter{algorithmctr}
\renewcommand{\thealgorithmctr}{\arabic{algorithmctr}}
   {\mbox{}\\*[\parskip]\begin{minipage}{\linewidth}%
       \refstepcounter{algorithmctr}\begin{list}{}{%
       \setlength{\rightmargin}{0\linewidth}%
       \setlength{\leftmargin}{.05\linewidth}}%
       \rmfamily\small
       \item[]{\setlength{\parskip}{0ex}\hrulefill\par%
        \nopagebreak{\bfseries\textsf{Algorithm \thealgorithmctr~}}}}%
   {{\setlength{\parskip}{-1ex}\nopagebreak\par\hrulefill\\*[2ex]\par}%
   \end{list}\end{minipage}}